\documentclass[11pt]{article}
\usepackage[T1]{fontenc}
\usepackage[margin=2.6cm]{geometry}
\usepackage{lmodern,microtype}
\usepackage{amsmath,amssymb,amsthm,mathtools}
\usepackage{enumitem,booktabs}
\usepackage[hidelinks]{hyperref}
\usepackage{cite}
\usepackage{tikz-cd}
\newtheorem{theorem}{Theorem}[section]

\newtheorem{proposition}[theorem]{Proposition}
\newtheorem{lemma}[theorem]{Lemma}

\newtheorem*{remark}{Remark}
\theoremstyle{definition}
\newtheorem{definition}[theorem]{Definition}

\newcommand{\G}{\mathcal G_p}
\newcommand{\States}{\mathcal A_p}
\newcommand{\sig}{\sigma}
\newcommand{\Types}{\mathcal X}
\newcommand{\Wild}{\operatorname{Wild}}
\newcommand{\wt}[1]{\widehat{#1}}

\setlist[enumerate]{label=(\roman*),itemsep=3pt}
\setlist[itemize]{itemsep=3pt}
\allowdisplaybreaks[1]
\title{Avoiding squares over words with lists of size three}
\author{Matthieu Rosenfeld\footnote{LIRMM, Université de Montpellier, CNRS, Montpellier, France}
\footnote{Laboratoire GAATI, University of French Polynesia}}

\begin{document}
\maketitle
\vspace{0.5em}
\begin{abstract}
We prove that every sequence of lists of size three admits an infinite square-free word, establishing that the Thue choice number of the path is three. More precisely, the number of compatible square-free words of length $n$ is at least $c(5/4)^n$, where $c>0$ is an absolute constant.

The proof combines a computer-assisted construction of weights with the weighted counting method previously introduced by the author.
We first consider words avoiding only squares of bounded period and use a finite state model to certify a uniform lower bound on their growth. 
The counting argument then controls the loss caused by excluding longer squares. The main new ingredient is a coarser description of the states, which makes the computation feasible without restricting the union of the lists.

The same computation yields a deterministic polynomial-time strategy for constructing an infinite square-free word when each list is chosen adversarially after the preceding letters have been fixed.
\end{abstract}

\section{Introduction}
A \emph{square} is a word $xx$ with $x$ nonempty. A word is \emph{square-free} if none of its factors is a square. 
For instance, $hotshots$ is a square while $minimize$ is square-free.
In 1906, Thue constructed an infinite ternary word avoiding squares \cite{Thue06,Thue1}.
Since then, many generalizations of Thue's result have been studied, including the list version of the problem.
Given a sequence of lists $L=(L_i)_{i\ge1}$, a word $w$ \emph{respects} $L$ if $w_i\in L_i$ for every position $i$ of $w$. 
The \emph{Thue choice number of the path},
denoted by $\pi_{\mathrm{ch}}(P)$, is the least integer $k$ such that for every infinite sequence $L$ of lists of size at least $k$, there exists an infinite square-free word respecting $L$.
The historical motivation for the study of the Thue choice number of the path is the generalization of the notion of nonrepetitive sequences to nonrepetitive colorings of graphs, which has received a lot of attention over the years (see for instance the recent survey of Wood \cite{woodSurvey}, or some of the results on the nonrepetitive chromatic number of graphs of bounded degree~\cite{Dujmovic2016,MontassierEntropie,HARANT2012374,rosenfeldCounting}).

Grytczuk, Przyby{\l}o and Zhu proved that $\pi_{\mathrm{ch}}(P)\le4$, and conjectured that $\pi_{\mathrm{ch}}(P)=3$ \cite{GrytczukLLLlist}.
This question, first mentioned in \cite{Grytczukpathquestion}, has appeared in a large number of papers, most often in its graph-theoretic form \cite{questionpathpathwidth,Grytczukpathquestion,Grytczukpathquestion2,GrytczukGame,rosenfeldCounting,Skrabuvlakova2015TheTC,woodSurvey,deuxplusepsilon,Bosek2021,KozikMicek2013,Przybylo2014,Peterin2018,Dujmovic2016,Fiorenzi2011,AlonSpencer2016,Cranston2024, R, rosenfeld2020far,cubeFreeWithVuong}.
In \cite{R}, I proved that lists of size three suffice when all lists are subsets of a common set of size $4$.
Here, we remove the restriction on the alphabet and settle the full conjecture.

\begin{theorem}\label{thm:main}
There is a constant $c>0$ such that, for every sequence $L=(L_i)_{i\ge1}$ of lists of size at least three and every $n\ge0$, there are at least
\[
 c\left(\frac{5}{4}\right)^n
\]
square-free words of length $n$ respecting $L$.

In particular, there exists an infinite square-free word respecting $L$, hence $\pi_{\mathrm{ch}}(P)=3$.
\end{theorem}

Our proof technique is essentially the one used in \cite{R}. 
In the conclusion of \cite{R}, I argued that this technique would probably settle the conjecture if one could carry out computations roughly $10^{7}$ times more expensive than those performed there\footnote{Under the stronger hypothesis that the growth rate is large enough. We are proving in the current article that the growth rate is large enough, so the previous proof from \cite{R} would also work, given enough computational power.}.
The main idea of that proof was to start by considering the language $\mathcal{L}_p$ of words over $\{1,2,3,4\}$ that avoid squares of length at most $2p$.
This language is defined by a finite set of forbidden factors, and can thus be recognized by a finite automaton, which can be used to compute a lower bound on the growth rate of $\mathcal{L}_p$.
More precisely, we computed a number $\alpha$, an analogue of the growth rate of $\mathcal{L}_p$ when an adversary can forbid one arbitrary transition out of each state (which corresponds to having only $3$ out of the $4$ letters available).
From there, we used the counting argument from \cite{rosenfeldCounting} and techniques from \cite{shurgrowthrate,Kolpakov2007} to bound the influence of longer squares.
For $p=21$, we could verify that $\alpha$ was large enough to satisfy our conditions.
When working with arbitrary alphabets, the only obstacle preventing this technique from proving the general conjecture was the computational cost of the automaton.

The technique from \cite{R} was also recently used in \cite{cubeFreeWithVuong} to prove that every sequence of lists of size $2$ admits an infinite cube-free word, with no restriction on the union of the lists.
The proof in \cite{cubeFreeWithVuong} is a direct adaptation of the proof from \cite{R} to the cube-free setting, but the computer verification is more manageable in this case, so no restriction on the union of the alphabets was needed.
Moreover, Bui and I gave a second proof that relies on exactly the same computer verification, but has slightly different consequences.
In particular, it implies that if the sequence of lists is computable, then there exists a computable cube-free word that respects it.
In fact, consider the setting where the sequence of lists is replaced by an adversary who chooses each list only after we have chosen the previous letters of the word. In this setting, there is  a strategy, computable in polynomial time\footnote{Here, we assume that comparing the letters at two positions can be done in polynomial time.
This in fact depends on the encoding of the letters, and we also need to impose that the adversary cannot play a completely arbitrary letter. 
Indeed, suppose that all the alphabets are subsets of $\mathbb{N}$ and that we use the standard binary encoding of the integers.
If the adversary provides the letter $2^{2^n}$ at step~$n$, then reading this letter already takes time exponential in~$n$.
This problem can be avoided by imposing that, at step~$n$, the
adversary provides letters bounded by~$3n$ (this does not change the combinatorial game, since we need to examine at most three letters per list and may rename letters in the order in which they first appear).
We will ignore this question in the rest of the article, and simply assume that, at step~$n$, any two letters of the word or of the next list can be compared in time polynomial in~$n$.} that produces a cube-free word respecting the lists.
This second technique also works here, and we get the following result.

\begin{theorem}\label{thm:onlineLists}
Consider the following game.
At each step $n+1$, an adversary reveals a set $L_{n+1}$ of at least three letters, which may depend on the letters $w_1,\ldots,w_n$ already chosen, and we choose $w_{n+1}\in L_{n+1}$.
There is a deterministic strategy, computing each move in time polynomial in $n$, that guarantees that the infinite word $w_1w_2\cdots$ is square-free.
\end{theorem}

In fact, this theorem is a corollary of the following result.
\begin{theorem}\label{thm:computationalGame}
There is a map $b$, computable in polynomial time, from finite words to letters, such that every infinite word $w_1w_2\cdots$ satisfying $w_{n+1}\notin\{w_n,\,b(w_1\cdots w_n)\}$ for all $n\ge1$ is square-free. The first letter $w_1$ can be chosen arbitrarily.
\end{theorem}
One way to understand this statement is through another game: at each step $n+1$ with $n\ge1$, given the current word $w_1\cdots w_n$, we can forbid two letters $x$ and $y$, and the adversary chooses $w_{n+1}$ with the only constraint that it differs from $x$ and $y$.
Theorem~\ref{thm:computationalGame} states that the map $b$ provides a strategy to keep the word square-free forever against any strategy from an opponent, when we start from the empty word.
Note that we clearly have to forbid $w_n$ at step $n+1$, since otherwise the adversary could choose $w_{n+1}=w_n$ and create a square, so the only real choice is the second forbidden letter, which is given by $b(w_1\cdots w_n)$.
This does not mean that we can win the game starting from any square-free word (after $0102010$, we would need to forbid $0$, $1$ and $2$ for the next letter). 
Forbidding the letter given by $b$ not only helps avoid squares at the next step, but also makes it easier to avoid squares in the future. 

The article is organized as follows. In Section~\ref{sec:count}, we present a weighted counting criterion that reduces the proof of Theorem~\ref{thm:main} to the existence of a function $f$ satisfying some local inequalities.
In Section~\ref{sec:otherconsequences}, we prove that the existence of the same function $f$ (with some extra computability conditions) also implies Theorem~\ref{thm:computationalGame}. The rest of the article is then devoted to the construction of such a function $f$.
In Section~\ref{sec:states}, we associate with every word a state, from a finite set that encode sufficient information to recover information about short squares. We then prove in Section~\ref{sec:types} and Section~\ref{sec:certificate}, that if we can find a function $f_s$ over the states that satisfies some local inequalities, then we can construct the desired function $f$ over the words. Under the assumption that this function $f_s$ exists, we can then deduce Theorems~\ref{thm:main}, \ref{thm:onlineLists} and~\ref{thm:computationalGame}.
Finally, in Section~\ref{sec:computations}, we provide some details about the implementation of the computations that verify the existence of the function $f_s$ for $p=21$. We conclude with some open problems in Section~\ref{sec:open}.

\section{A weighted counting argument}\label{sec:count}

Fix $p\ge2$, and let $\G$ be the set of finite words containing no
square of period at most $p$, where the period of $xx$ is $|x|$.
The following form of the counting argument isolates the properties
needed from the weights.

\begin{theorem}\label{thm:criterion}
Suppose that a function $f:\G\to[0,M]$, with $M<\infty$ and
$f(\varepsilon)>0$, satisfies the following conditions for some
$\alpha>0$:
\begin{enumerate}
\item\label{item:local} For every $w\in\G$ and every set $B$ of three letters,
\begin{equation}\label{eq:local}
 \sum_{\substack{a\in B\\wa\in\G}} f(wa)\ge\alpha f(w).
\end{equation}
\item\label{item:delete} If $yxx\in\G$ and $|x|>p$, then $f(yxx)=f(yx)$.
\end{enumerate}
If $\beta>1$ satisfies
\begin{equation}\label{eq:criterion}
 \alpha-\frac{\beta^{1-p}}{\beta-1}\ge\beta,
\end{equation}
then every list assignment with lists of size three has at least
$f(\varepsilon)\beta^n/M$ compatible square-free words of length $n$, for every $n\ge0$.
\end{theorem}

\begin{proof}
Fix the lists, and let $S_n$ be the set of compatible square-free words of length $n$.
For a finite set $U\subseteq\G$, write $\wt U=\sum_{w\in U}f(w)$. We prove by induction that for all $i\ge0$,
\begin{equation}\label{eq:induction}
 \wt S_{i+1}\ge\beta\wt S_i.
\end{equation}
Suppose that \eqref{eq:induction} holds for all $i<n$. We will prove that \eqref{eq:induction} holds for $i=n$.

Let $H_{n+1}$ consist of the words obtained by extending a word of $S_n$ by a letter of $L_{n+1}$ without creating a square of period at most $p$.
Then by \eqref{eq:local}, we get
\begin{align*}
\wt H_{n+1}
&=\sum_{w\in S_n}\sum_{\substack{a\in L_{n+1}\\wa\in\G}}f(wa)
 \ge\sum_{w\in S_n}\alpha f(w)=\alpha\wt S_n.
\end{align*}
Some of these extensions do not bellong to $S_{n+1}$, if they create a square of period larger than $p$. 
Precisely, every word in $H_{n+1}\setminus S_{n+1}$ has a square suffix of period larger than $p$, since its prefix of length $n$ is square-free.
For $p<i\le\lfloor(n+1)/2\rfloor$, let $F_i$ be the set of words in $H_{n+1}$ ending in a square of period $i$.
Then, $S_{n+1}=H_{n+1}\setminus\bigcup_{i=p+1}^{\lfloor(n+1)/2\rfloor}F_i$, and
\begin{equation}\label{eq:subtract}
 \wt S_{n+1}\ge\alpha\wt S_n-
 \sum_{i=p+1}^{\lfloor(n+1)/2\rfloor}\wt F_i.
\end{equation}

By definition, every $w\in F_i$ can be written as $w=yxx$ with $|x|=i>p$, and since every proper prefix of $w$ is square-free, we have $yx\in S_{n+1-i}$. Thus the map
\[
F_i\longrightarrow S_{n+1-i},\qquad yxx\longmapsto yx
\]
is well-defined and injective, since $w$ can be recovered from $yx$ by appending its suffix of length $i$. Moreover, condition~\ref{item:delete} gives $f(yxx)=f(yx)$, so this map
preserves weights. This implies,
\[
\wt F_i\le \wt S_{n+1-i}\le\beta^{1-i}\wt S_n,
\]
where the last inequality follows from applying the induction hypothesis \eqref{eq:induction} for all $j\in \{n+1-i,\ldots, n-1\}$.
Substituting in~\eqref{eq:subtract},
we obtain
\[
 \wt S_{n+1}
 \ge \alpha\wt S_n-\sum_{i=p+1}^{\lfloor(n+1)/2\rfloor}\beta^{1-i}\wt S_n
 \ge \left(\alpha-\sum_{i=p+1}^{\infty}\beta^{1-i}\right)\wt S_n
 =\left(\alpha-\frac{\beta^{1-p}}{\beta-1}\right)\wt S_n
 \ge\beta\wt S_n\,
\]
where the last inequality is simply our hypothesis \eqref{eq:criterion}.
This concludes the proof that for all $i\ge0$, $\wt S_{i+1}\ge\beta\wt S_i$.
Finally, $M|S_n|\ge\wt S_n\ge f(\varepsilon)\beta^n$, as desired.
\end{proof}

The function $f$ might at first seem a bit abstract, and it is not clear how to prove the existence of such a function. However, avoiding squares of period at most $p$ is a local property, and the function $f$ can be constructed so that it only depends on the suffix of length at most $2p$ of its input. Moreover, by symmetry all letters play a similar role, hence it is natural to ask that the function $f$ only depend on the last $2p$ letters, up to renaming. For a given $p$, there are only finitely many possibilities for the last $2p$ letters of a word, up to renaming. The best way to construct $f$ is to construct a set of states that encode enough information about the last $2p$ letters of a word, up to renaming, and then define $f$ as a function on these states. This was already the approach used in \cite{R}, and Theorem \ref{thm:criterion} is simply a more abstract rewriting of \cite[Theorem 5]{R}.

In \cite{R}, the states were ``exact'' in the sense that they encoded enough information so that there was a bijection between the transitions of a word and the transitions of its state. The set of states was too large to fit in a computer, and we had to settle for a weaker result. In the current paper, we will use a coarser set of states, which loses some of this information for the states that have many outgoing transitions. This implies that the constant $\alpha$ will be slightly smaller than with exact states, but this effect is negligible while making the set of states much smaller and allowing us to finish the computations.

In the remainder of the proof we construct $f$ as a function of a finite ``state'' depending on the last $2p$ letters of the word. This is done in Sections~\ref{sec:states} to~\ref{sec:computations}. Before that, we prove in Section~\ref{sec:otherconsequences} that the existence of such a function $f$ also implies Theorem~\ref{thm:computationalGame}.

\section{The computational consequences of the same condition}\label{sec:otherconsequences}
For the algorithmic statements, we take the alphabet to be $\mathbb N$.
Given a function $f:\G\to[0,M]$ and a real number $\beta>1$, we define the potential $\omega_{f,\beta}(u)$ of any word $u\in \G$ as follows:
\[
    \omega_{f,\beta}(u)=
    \sum_{\substack{v\text{ such that: }\\uv\in\G \\ uv=xyy \text{ with }|v|<|y|}} f(uv)\beta^{-|v|}\,.
\]
This function is the main ingredient used to define the function $b$ from Theorem~\ref{thm:computationalGame}.
First note that, if $u$ ends with a square, then one of the summands corresponds to $v=\varepsilon$, and thus $\omega_{f,\beta}(u)\ge f(u)$. This implies the following lemma.
\begin{lemma}\label{lem:oneStepSquareFree}
If $u$ is square-free and $a$ is a letter such that $ua\in\G$ and $\omega_{f,\beta}(ua)< f(ua)$, then $ua$ is square-free.
\end{lemma}

The core of our proof is that under our assumptions, if $\omega_{f,\beta}(u)< f(u)$ then at most two letters $a$ satisfy $ua\notin\G$ or $\omega_{f,\beta}(ua)\ge f(ua)$. By forbidding these two letters, we can ensure that $ua\in\G$ and $\omega_{f,\beta}(ua)<f(ua)$ at the next step, and we can then apply the same argument inductively.
We now prove one step of the induction in the following lemma.

\begin{lemma}\label{lem:computabilityOneStep}
Let $f$ and $\beta$ satisfy the hypotheses of
Theorem~\ref{thm:criterion}, and let $u\in\G$ be such that
$\omega_{f,\beta}(u)<f(u)$. Then every set of three letters contains a
letter $a$ such that $ua\in\G$ and $\omega_{f,\beta}(ua)<f(ua)$.
Equivalently, at most two letters $a$ satisfy $ua\notin\G$ or
$\omega_{f,\beta}(ua)\ge f(ua)$.
\end{lemma}
\begin{proof}
We first bound the sum of the potentials of all the extensions of $u$ by a
letter\footnote{Although this sum might seem infinite,
$\omega_{f,\beta}(ua)=0$ whenever $a$ does not occur in $u$, so it has
finitely many nonzero terms (see the proof of Lemma~\ref{lem:computableStrategy}).}:
\begin{align*}
\sum_{\substack{a\\ua\in\G}} \omega_{f,\beta}(ua)
&=\sum_{\substack{a\\ua\in\G}}\sum_{\substack{v\text{ such that: }\\uav\in\G \\ uav=xyy \text{ with } |v|<|y|}} f(uav)\beta^{-|v|}\\
&\le\sum_{\substack{v'\text{ such that: }\\uv'\in\G\\ uv'=xyy \text{ with }1\le|y|, |v'|\le|y|}} f(uv')\beta^{1-|v'|}\\
&\le \beta \omega_{f,\beta}(u)+\sum_{\substack{v'\text{ such that: }\\uv'\in\G\\ uv'=xyy \text{ with }1\le|y|, |v'|=|y|}} f(uv')\beta^{1-|v'|}\,,
\end{align*}
where the last inequality is obtained by splitting the sum according to whether $|v'|<|y|$ or $|v'|=|y|$, and noticing that the terms with $|v'|<|y|$ form a subset of the terms in the definition of $\omega_{f,\beta}$.
In the last sum, we have $uv'=xyy$ with $|v'|=|y|$, so $u=xy$ and $v'=y$. Since $uv'\in\G$, we have $|y|>p$, so condition~\ref{item:delete} gives $f(uv')=f(u)$. Moreover, for each positive integer $\ell$ there exists at most one such $v'$ of length $\ell$, namely the suffix of length $\ell$ of $u$. Thus we get
\begin{align*}
\sum_{\substack{a\\ua\in\G}} \omega_{f,\beta}(ua)
&\le \beta \omega_{f,\beta}(u)+\sum_{\ell = p+1}^{\infty} f(u)\beta^{1-\ell}\\
&= \beta \omega_{f,\beta}(u)+\frac{\beta^{1-p}}{\beta-1}f(u)\\
&\le \beta \omega_{f,\beta}(u)+f(u)(\alpha-\beta)\,,
\end{align*}
where the last inequality follows from condition~\eqref{eq:criterion}. 
Now let $B$ be a set of three letters. Since $\omega_{f,\beta}(u)<f(u)$, the preceding bound and~\eqref{eq:local} give
\[
\sum_{\substack{a\in B\\ua\in\G}}\omega_{f,\beta}(ua)
\le\sum_{\substack{a\\ua\in\G}}\omega_{f,\beta}(ua)
<\alpha f(u)
\le\sum_{\substack{a\in B\\ua\in\G}}f(ua)\,.
\]
This implies that there is at least one letter $a\in B$ such that $ua\in\G$ and $\omega_{f,\beta}(ua)<f(ua)$, as desired. Hence at most two letters $a$ satisfy $ua\notin\G$ or $\omega_{f,\beta}(ua)\ge f(ua)$, as we could have chosen $B$ to be the set of these three letters otherwise.
\end{proof}
We now prove Theorem~\ref{thm:computationalGame} as a corollary of Lemma~\ref{lem:computabilityOneStep}, under the assumption that there exist $f$ and $\beta$ satisfying the hypotheses of Theorem~\ref{thm:criterion}, with $f$
computable in polynomial time and independent of the names of the letters. We check in Section~\ref{sec:certificate} that the function $f$ constructed there has these properties.

\begin{lemma}\label{lem:computableStrategy}
Let $f$ and $\beta$ satisfy the hypotheses of Theorem~\ref{thm:criterion}.
Assume moreover that $\beta$ is rational, that $f$ takes rational values
and is computable in polynomial time, and that $f$ does not depend on the
names of the letters: $f(\pi(w))=f(w)$ for every $w\in\G$ and every
injective map $\pi$ from letters to letters, applied letter by letter.
Then there is a map $b$, computable in polynomial time, from finite
words to letters, such that every infinite word
$w_1w_2\cdots$ satisfying $w_{n+1}\notin\{w_n,\,b(w_1\cdots w_n)\}$ for
all $n\ge1$ is square-free. The first letter $w_1$ can be chosen arbitrarily.
\end{lemma}
\begin{proof}
For $u\in\G$, let $\mathrm{Bad}(u)$ be the set of letters $a$ such that
$ua\notin\G$ or $\omega_{f,\beta}(ua)\ge f(ua)$. 

For every nonempty $u\in\G$ with $\omega_{f,\beta}(u)<f(u)$, we have $|\mathrm{Bad}(u)|\le2$ by Lemma~\ref{lem:computabilityOneStep}. 
Moreover, since $\omega_{f,\beta}(u)<f(u)$, every letter not occurring in $u$ lies outside $\mathrm{Bad}(u)$. Indeed, for such a letter $a$, we have $ua\in\G$ and $\omega_{f,\beta}(ua)=0$.
If $a$ is the last letter of $u$, then $a\in\mathrm{Bad}(u)$, since $ua\notin\G$. Hence $|\mathrm{Bad}(u)\setminus\{a\}|\le1$, and we define $b(u)$ to be this other bad letter if it exists, and $0$ otherwise\footnote{Note that the choice of $0$ is purely arbitrary here, we simply need to fix a convention for the sake of computability. The same comment is true for the next sentence.}. For every other finite word $u$, set $b(u)=0$.

Let $w_1w_2\cdots$ satisfy $w_{n+1}\notin\{w_n,b(w_1\cdots w_n)\}$ for all $n$. 

We prove by induction that every prefix $u$ is square-free and satisfies $\omega_{f,\beta}(u)<f(u)$. This holds for $u=\varepsilon$, since $\omega_{f,\beta}(\varepsilon)=0<f(\varepsilon)$. 
If it holds for $u=w_1\cdots w_n$, then $w_{n+1}\notin\mathrm{Bad}(u)$, so $uw_{n+1}\in\G$ and $\omega_{f,\beta}(uw_{n+1})<f(uw_{n+1})$, and $uw_{n+1}$ is square-free by Lemma~\ref{lem:oneStepSquareFree}.

It remains to verify that $b$ is computable in polynomial time. Let $u$ be a finite word of length $n$. If $u=\varepsilon$ or $u\notin\G$, return $0$.
Otherwise, in the definition of $\omega_{f,\beta}(u)$, the word $v$ is determined by $|y|$ and $|v|$. Enumerating these pairs and discarding duplicate completions $v$ gives a sum of $O(n^2)$ terms, each computable in polynomial time.
Since $\beta$ and the values of $f$ are rational, these terms are rational numbers of polynomial size, and $\omega_{f,\beta}(u)$ can be compared exactly with $f(u)$ in polynomial
time. 
If $\omega_{f,\beta}(u)\ge f(u)$, return $0$. Otherwise, every letter of $\mathrm{Bad}(u)$ already occurs in $u$. For each such occurring letter $a$, first test whether $ua\in\G$ and, if so, compare $\omega_{f,\beta}(ua)$ with $f(ua)$. This determines $\mathrm{Bad}(u)$ and $b(u)$ in polynomial time.
\end{proof}

\section{Finite states and deletion of long squares}\label{sec:states}
For a word $w=w_1\cdots w_n$ with $n\ge p$, write $w[j]=w_{n+1-j}$
for its letter at distance $j$ from the right end (in particular, $w[1]$ is the last letter of $w$). For $1\le q\le p$, let
\[
 r_q(w)=\max\bigl\{r\in\{0,\ldots,n-q\}:
          w[j]=w[j+q]\text{ for }1\le j\le r\bigr\}.
\]
Thus $r_q(w)$ is the largest integer $r$ such that the suffix of length $q+r$ of $w$ has period $q$. If $w\in\G$, then $r_q(w)\le q-1$.
From this definition, for any letter $a$,
\begin{equation}\label{eq:overlap-update}
 r_q(wa)=
 \begin{cases}
 r_q(w)+1,&q\le n\text{ and }a=w[q],\\
 0,&\text{otherwise}.
 \end{cases}
\end{equation}
In particular, for $w\in\G$, and any letter $a$,
\begin{equation}\label{eq:legal-word}
 wa\in\G
 \quad\Longleftrightarrow\quad
 \forall q\le\min(n,p)\text{ such that }a=w[q],\quad r_q(w)<q-1.
\end{equation}
Indeed, a square of period at most $p$ in $wa$ is necessarily a suffix of $wa$, and if it has period $q$ then $a=w[q]$ and $r_q(w)\ge q-1$. 

To avoid treating short words separately, fix $p+1$ letters $c_0,c_1,\ldots,c_p$ that do not belong to the alphabet $\Sigma$ of the lists, and put $\rho=c_p\cdots c_1c_0$. 
Each symbol of $\rho$ occurs exactly once in $\rho w$ and never in $w$, so no square of $\rho w$ contains a symbol of $\rho$.
Hence $\rho w\in\G$ if and only if $w\in\G$, and the squares of $\rho w$ are those of $w$. We will only use the states of the words $\rho w$, which have length at least $p+1$.

All letters play the same role, so we can always assume up to symmetry that the last two letters of a word are $1$ and $0$. 
In order to decrease the number of states, we will remember the positions of all occurrences of $0$ and $1$, but we will forget the details about all other letters. 
To avoid losing too much information, we will still keep track of the overlaps of the last $p$ letters given by the different values of $r_q$ for $1\le q\le p$. 
More precisely, for $w\in\G$ with $|w|\ge p$ and $1\le j\le p$, let
\[
 \lambda_j(w)=
 \begin{cases}
 0,&w[j]=w[1],\\
 1,&w[j]=w[2],\\
 \star,&\text{otherwise}.
 \end{cases}
\]
Define the \emph{state} of $w$ by
\begin{equation}\label{eq:state}
 \sig(w)=\bigl(\lambda_p(w),\ldots,\lambda_1(w);
                         r_p(w),\ldots,r_1(w)\bigr).
\end{equation}
The label $\star$ does not mean that all the corresponding positions
carry the same letter. 
It means that the only information that we have to distinguish these letters consists of the values of the $r_q(w)$: the maximality condition in the definition of $r_q(w)$ implies that the letters at positions $r_q(w)+1$ and $r_q(w)+1+q$ are different, whenever $r_q(w)+1+q\le n$.

For definiteness, let $\States$ be the finite set of tuples $A=(\lambda_p,\ldots,\lambda_1;r_p,\ldots,r_1)$  with $0\le r_q\le q-1$, labels in $\{0,1,\star\}$, $\lambda_1=0$ and $\lambda_2=1$. 
Every word in $\G$ of length at least $p$ has a state in $\States$, although not every tuple in $\States$ is realized by a word.
The state of the prefix $\rho$, which we call the \emph{start state}, is
\[
 A_0=\sig(\rho)=(\star,\ldots,\star,1,0;0,\ldots,0).
\]

To define the function $f$ from Theorem \ref{thm:criterion}, we will define a weight function $D:\States\to\mathbb Z_{\ge0}$, and then set $f(w)=D(\sig(\rho w))$. 
The following lemma ensures that any function $f$ defined in this way has the second property needed to apply Theorem \ref{thm:criterion}, that is, if $yxx\in\G$ and $|x|>p$, then $f(yxx)=f(yx)$.
\begin{lemma}\label{lem:deletion}
If $yxx\in\G$ and $|x|>p$, then $\sig(yxx)=\sig(yx)$.
\end{lemma}
\begin{proof}
Write $u=yxx$, $v=yx$ and $i=|x|$. We first claim that for all $q\in\{1,\ldots, p\}$,
\begin{equation}\label{eq:overlap-bound}
 q+r_q(u)< |x|.
\end{equation}
Otherwise, $x$ has period $q$.
Put $d=i-q$, so $1\le d\le q-1$. The prefix of $x$ of length $d$ equals its suffix of length $d$. 
These two copies meet at the junction in $xx$, forming a square of period $d<p$, contrary to $u\in\G$.

For all $q\in\{1,\ldots, p\}$,  $q+r_q(u)+1\le |x|$ implies that the suffix of length $q+r_q(u)$ of $x$ has period $q$, but not the suffix of length $q+r_q(u)+1$. 
So the value of $r_q(u)$ is uniquely determined by $x$, independently of $y$. In particular, $r_q(u)=r_q(v)$ for all $q\in\{1,\ldots, p\}$.
The labels of the states also agree, as they depend only on the suffix of length $p<|x|$.
\end{proof}

\section{Extensions and possible letter types}\label{sec:types}
We now define state transitions corresponding to the addition of a letter. 
For every letter $a$ and word $w\in\G$ with $|w|\ge p$, define the \emph{type} of $a$ after $w$
to be its occurrence set in the suffix of length $p$:
\[
 \mathcal{T}_w(a)=\{j\in\{1,\ldots,p\}:w[j]=a\}.
\]
Distinct letters have disjoint types. 
The type of $a$ after $w$ can be seen as choosing a letter from the suffix, or a letter absent from it. 
We will define a corresponding notion of type for states, which will be a set of positions in the suffix of length $p$ encoded by the state.

For a state $A=(\lambda_p,\ldots,\lambda_1;r_p,\ldots,r_1)$, we let
\[
 Z_0(A)=\{j:\lambda_j=0\},\qquad
 Z_1(A)=\{j:\lambda_j=1\},\qquad
 \Wild(A)=\{j:\lambda_j=\star\}.
\]
Consider a set $X$ equal to $Z_0(A)$ or
$Z_1(A)$, or contained in $\Wild(A)$. The type $X$ is \emph{allowed} if $r_q<q-1$ for every $q\in X$. In particular, by this definition the type $Z_0(A)$ is never allowed. Intuitively, a type is allowed if adding the corresponding letter does not create a square of period at most $p$ (see Lemma~\ref{lem:update}).
For an allowed $X$, we define the successor 
\[
\Phi(A,X)=(\lambda'_p,\ldots,\lambda'_1;r'_p,\ldots,r'_1)
\]
where
\begin{equation}\label{eq:state-overlaps}
 r'_q=\begin{cases}r_q+1,&q\in X,\\0,&q\notin X,\end{cases}
 \qquad 1\le q\le p,
\end{equation}and, for $1\le j\le p$,
\begin{equation}\label{eq:state-labels}
 \lambda'_j=
 \begin{cases}
 0,&j=1\text{ or }j-1\in X,\\
 1,&\lambda_{j-1}=0,\\
 \star,&\text{otherwise}.
 \end{cases}
\end{equation}

The following lemma shows that the successor of a state is well-defined, and that it corresponds to the state of a word after adding a letter of the corresponding type.

\begin{lemma}\label{lem:update}
Let $w\in\G$ be a word of length at least $p$, and $a$ be a letter.
Then $wa\in\G$ if and only if $\mathcal{T}_w(a)$ is allowed at $\sig(w)$.
When it is allowed,
\[
 \sig(wa)=\Phi(\sig(w),\mathcal{T}_w(a)).
\]
\end{lemma}

\begin{proof}
The set $\mathcal{T}_w(a)$ records exactly the positions occupied by $a$ among the last $p$ letters. 
Legality and the overlap update therefore follow from~\eqref{eq:legal-word} and~\eqref{eq:overlap-update}.
After a legal extension, the distinguished letters are $a$ and the previous last letter. 
The positions in $\mathcal{T}_w(a)$ acquire label $0$, the former $0$-positions acquire label $1$, and all other positions acquire label $\star$. 
Shifting the positions by one gives~\eqref{eq:state-labels}.
\end{proof}

We have proven that the states behave nicely, but for any given state we need to be able to tell what type of letter is possible after the state.
A letter that is not one of the last two letters of $w$ is called an \emph{outside letter}: its type is contained in $\Wild(\sig(w))$, and we call the subsets of $\Wild(A)$ \emph{outside types}. However, some subset of $\Wild(A)$ do not correspond to the occurrence set of an outside letter. The labels of a state only contain information about the positions of $0$ and $1$ (the two types $Z_0(A)$ and $Z_1(A)$ are easily defined). However, we still have some partial information about the other positions which impose some pair of letters to be identical or different.

The overlaps given by the $(r_q)_{1\le q\le p}$ force certain $\star$ positions to carry the same letter. Similarly, choosing an outside type $X\subseteq\Wild(A)$ forces all positions of $X$ to carry the same letter.
If all of these constraints impose the presence of a square, then $X$ cannot be the occurrence set of an outside letter (since the suffix of length $p$ of a word in $\G$ is square-free).
Moreover, the $r_q$ are maximal: whenever $r_q+1+q\le p$, the letters at positions $r_q+1$ and $r_q+1+q$ are distinct, so an outside type cannot include both of these positions.

Formally, for $X\subseteq\Wild(A)$, let $\equiv_{A,X}$ be the equivalence relation on $\{1,\ldots,p\}$ generated by
\begin{equation}\label{eq:forced-equality}
\begin{aligned}
 &j\equiv_{A,X}j+q &&(1\le q\le p,\ 1\le j\le r_q,\ j+q\le p),\\
 &j\equiv_{A,X}j' &&(\lambda_j=\lambda_{j'}\in\{0,1\}\text{, or }j,j'\in X).
\end{aligned}
\end{equation}
That is, two positions are equivalent for $\equiv_{A,X}$ if they are forced to carry the same letter by the overlaps $r_q$, by the labels, or because they both belong to $X$. Let
$s_X=s_X(p)\cdots s_X(1)$ be the word whose letter $s_X(j)$ is the
$\equiv_{A,X}$-class of~$j$, and define
\begin{equation}\label{eq:outside-types}
 \Types_{\mathrm{out}}(A)\coloneq
 \left\{X\subseteq\Wild(A):
 \begin{array}{l}
 X\text{ is a union of }{\equiv_{A,\emptyset}}\text{-classes},\\
 s_X\text{ is square-free, and}\\
 r_q+1\not\equiv_{A,X}r_q+1+q\text{ whenever }r_q+1+q\le p
 \end{array}\right\},
\end{equation}
and
\[
 \Types(A)\coloneq\{Z_0(A),Z_1(A)\}\cup\Types_{\mathrm{out}}(A).
\]

 We now prove that the set of types of a word is contained in the set of types of its state.
\begin{lemma}\label{lem:outside}
If $w\in\G$ with $|w|\ge p$ and $a$ is a letter, then
$\mathcal{T}_w(a)\in\Types(\sig(w))$.
\end{lemma}
\begin{proof}
If $a$ is one of the last two letters of $w$, then $\mathcal{T}_w(a)=Z_0(\sig(w))$ or $\mathcal{T}_w(a)=Z_1(\sig(w))$, so the assertion is immediate. Otherwise, $\mathcal{T}_w(a)\subseteq\Wild(\sig(w))$. The two positions of each generating pair of~\eqref{eq:forced-equality} carry the same letter in~$w$:
by definition of $r_q$, by definition of the labels, or because both carry~$a$. 
Hence any two $\equiv_{\sig(w),\mathcal{T}_w(a)}$-equivalent positions carry the same letter. 
In particular, a position equivalent to a position of~$\mathcal{T}_w(a)$ carries~$a$, so it lies in~$\mathcal{T}_w(a)$, and $\mathcal{T}_w(a)$ is a union of $\equiv_{\sig(w),\emptyset}$-classes. 
A square in $s_{\mathcal{T}_w(a)}$ would give a square of period at most $p$ in~$w$, which is impossible since $w\in\G$. 
Finally, if $r_q+1+q\le p$, the maximality of $r_q$ gives $w[r_q+1]\neq w[r_q+1+q]$, so these two positions are not equivalent.
\end{proof}
\begin{remark}\label{rem:longer-window}
The relation~\eqref{eq:forced-equality} and the word $s_X$ only use the last $p$ positions, but the $r_q$ provide information about the positions up to $2p-1$. 
Keeping $X\subseteq\Wild(A)$ and the first condition in~\eqref{eq:outside-types} unchanged, one may extend $\equiv_{A,X}$ and $s_X$ to the positions $1,\ldots,2p-1$, using all the pairs $j\equiv_{A,X}j+q$ with $1\le j\le r_q$. 
Requiring the extended word to be square-free eliminates a few more types.
For instance, let $p=6$ and $w=cadaedba$. 
Then $\sig(w)=(\star,0,\star,\star,1,0;\;1,0,1,0,0,0).$ Definition~\eqref{eq:outside-types} gives
\[
 \Types_{\mathrm{out}}(\sig(w))=\bigl\{\emptyset,\{3\},\{4\},\{6\},\{3,6\},\{4,6\}\bigr\}.
\]
In particular, the type $X=\{4,6\}$ passes the test on the last six positions, where $s_X$ reads $\mathsf{x}\,0\,\mathsf{x}\,\kappa\,1\,0$, writing $\mathsf{x}$, $0$, $\kappa$ and $1$ for the classes of $X$, of the positions labelled $0$, of position~$3$ and of the positions labelled~$1$. 
But $r_6=1$ forces the letter at position~$7$ to be the last letter, so the positions $7,6,5,4$ read $0\,\mathsf{x}\,0\,\mathsf{x}$, a square of period~$2$. 
No word with this state has a letter occupying exactly the positions $4$ and~$6$, and the extended definition would remove this type.
 
Lemma~\ref{lem:outside} holds for this extended definition, with the same proof. 
On the one hand, this experimentally reduces the number of reachable states by about $13\%$ and improves the certified constant by less than $2\cdot10^{-5}$. On the other hand, each state costs about $20\%$ more to process and, since the tables are allocated by ``doubling'', this does not improve the memory or running time. We therefore stick to the simpler definition.
\end{remark}

\section{The finite weight inequalities}\label{sec:certificate}

Let $D:\States\to\mathbb Z_{\ge0}$. For a state $A$ and a type $t\in \Types(A)$, define
its extension weight by
\[
 v_D(A,t)=
 \begin{cases}
 D(\Phi(A,t)),&t\text{ is allowed},\\
 0,&\text{otherwise}.
 \end{cases}
\]
An \emph{admissible triple} at $A$ is a multiset $\{X_1,X_2,X_3\}$ of elements of $\Types(A)$ such that
$X_i\cap X_j=\emptyset$ whenever $i\ne j$.
Thus a nonempty set cannot occur twice, whereas $\emptyset$ may
occur several times. With the convention $\min\emptyset=+\infty$, let
\begin{equation}\label{eq:operator}
 T_D(A)=\min_{\text{admissible triples }\{t_1,t_2,t_3\}}
                    \sum_{j=1}^3 v_D(A,t_j).
\end{equation}

\begin{definition}\label{def:certificate}
A \emph{weight certificate} with constant $\alpha>0$ is a function
$D:\States\to\mathbb Z_{\ge0}$ such that
\begin{itemize}
\item $D(A_0)>0$, and
\item $T_D(A)\ge\alpha D(A)$ for every $A\in\States$.
\end{itemize}
\end{definition}

\begin{proposition}\label{prop:certificate}
If $D$ is a weight certificate with constant $\alpha$, then
$f(w)=D(\sig(\rho w))$ satisfies the hypotheses of
Theorem~\ref{thm:criterion}, with the same constant $\alpha$ and $M=\max_A D(A)$.
\end{proposition}
\begin{proof}
First note that $\States$ is finite, hence $D(A)$ and $f(w)$ are bounded. Moreover, $f(\varepsilon)=D(A_0)>0$. So we only need to check the two conditions (i) and (ii) of Theorem~\ref{thm:criterion}.

Let $w\in\G$, so that $\rho w\in\G$, and let $B\subseteq\Sigma$
have size three.
Since $B$ contains three distinct letters, the types of these letters after $\rho w$ are pairwise disjoint, and by Lemma~\ref{lem:outside} each type belongs to $\Types(\sig(\rho w))$. They therefore form an admissible triple at $\sig(\rho w)$.
Using Lemma~\ref{lem:update}, we obtain
\begin{align*}
 \sum_{\substack{a\in B\\wa\in\G}}f(wa)
 &=\sum_{\substack{a\in B\\\rho wa\in\G}}D(\sig(\rho wa))\\
 &=\sum_{\substack{a\in B\\\rho wa\in\G}}D(\Phi(\sig(\rho w),\mathcal{T}_{\rho w}(a)))\\
 &=\sum_{a\in B}v_D(\sig(\rho w),\mathcal{T}_{\rho w}(a))\\
 &\ge T_D(\sig(\rho w))\\
 &\ge\alpha D(\sig(\rho w))\\
 &=\alpha f(w).
\end{align*}
If $yxx\in\G$ is over $\Sigma$ and $|x|>p$, then $\rho yxx\in\G$,
so Lemma~\ref{lem:deletion} gives
\[
 f(yxx)=D(\sig(\rho yxx))=D(\sig(\rho yx))=f(yx),
\]
which concludes our proof.
\end{proof}

What remains is to find a weight certificate with a sufficiently large constant $\alpha$ to apply Theorem~\ref{thm:criterion}. 
\begin{proposition}[Computer-assisted]\label{prop:computation}
For $p=21$, there is a weight certificate with constant
\[
 \alpha=\frac{1300544}{1000000}=1.300544.
\]
\end{proposition}
This proposition is proved by computer, and we provide the details of the computation in the next section. First we show how to deduce Theorem~\ref{thm:main} and Theorem~\ref{thm:computationalGame} from it. 
\begin{proof}[Proof of Theorem~\ref{thm:main}]
First reduce every list to size three. Take the certificate of
Proposition~\ref{prop:computation} and set $\beta=5/4$.
Condition~\eqref{eq:criterion} is equivalent to
\[
 (\alpha-\beta)(\beta-1)\beta^{20}\ge1,
\]
which is easily verified to hold with the given values of $\alpha$ and $\beta$.
Proposition~\ref{prop:certificate} and Theorem~\ref{thm:criterion}
now give the claimed bound, with
$c=D(A_0)/\max_A D(A)>0$, independently of the lists.

The compatible square-free words form a finitely branching tree under extension by one letter, with a vertex at every depth. By K\"onig's lemma, this tree has an infinite branch, giving an infinite square-free word respecting the lists.
\end{proof}

\begin{proof}[Proof of Theorem~\ref{thm:computationalGame}]
Let $f(w)=D(\sig(\rho w))$, where $D$ is the certificate of Proposition~\ref{prop:computation}, and let $\beta=5/4$. As in the proof of Theorem~\ref{thm:main}, $f$ and $\beta$ satisfy the hypotheses of Theorem~\ref{thm:criterion}. Moreover, $\beta$ is rational and $D$ takes integer values. The state $\sig(\rho w)$ is determined by the last $2p$ letters of $\rho w$, since $r_q(\rho w)\le q-1$ for every $q\le p$, and $D$ is given by a finite table, so $f$ is computable in polynomial time. Finally, $\sig(\rho w)$ only depends on which positions of $\rho w$ carry the same letter, and the letters of $\rho$ do not belong to $\Sigma$, so $f$ does not depend on the names of the letters. Lemma~\ref{lem:computableStrategy} therefore gives Theorem~\ref{thm:computationalGame}.
\end{proof}

\section{Computation of the weight certificate}\label{sec:computations}
The proof of Proposition~\ref{prop:computation} is computer-assisted. The code is available at \url{https://github.com/matthieu-rosenfeld/Avoiding-squares-over-lists-of-size-3}.

The first critical part of the verification is to be able to generate the types and successors of a given state. To generate the types of a state $A$ we essentially use the definition of $\Types(A)$, which is a finite check. Practical details involve a few extra tricks like precomputing the equivalence classes of $\equiv_{A,\emptyset}$ and some information about the square-freeness to be able to efficiently generate and validate all the types of a state. The computation of the successor $\Phi(A,X)$ is straightforward using \eqref{eq:state-overlaps} and \eqref{eq:state-labels}. Crucially, we do not store the set of types associated with a state, but rather generate it each time we need it. This is much slower, but it saves a lot of memory and memory is the main limiting factor in our computation.

For the set of states, we do not consider the set of all states, but only the set $R$ of states reachable from the start state $A_0$, that is, obtained from $A_0$ by successive moves $A\mapsto\Phi(A,X)$ with $X\in\Types(A)$ allowed. For this we simply use a breadth-first search starting from $A_0$ and generate all the successors of each state. For all other states, we implicitly set $D(A)=0$, which is allowed by the definition of a weight certificate.

From there, we also need to be able to compute the operator $T_D(A)$ for each state $A$. Since there are finitely many types, this could be done by enumerating all triples of types. However, we can make this much more efficient by noting that the type $Z_0(A)$ is always forbidden and has zero extension weight, so it is always optimal to include it in a minimizing triple. So we only need to enumerate triples of types that contain $Z_0(A)$, which reduces the number of triples to check significantly.

In order to find the weights $D(A)$, we initially give each state in $R$ the integer weight $S=10^8$. At each iteration, we compute $T_D(A)$ for every $A\in R$. Unless the stopping criterion is met, we then normalize by setting
\[
 D_{\mathrm{new}}(A)=\left\lfloor\frac{S T_D(A)}{\max_{B\in R}T_D(B)}\right\rfloor.
\]
If the denominator is zero, the iteration stops.
Note that this is easily parallelizable since for the most part we read $D(A)$ and write $T_D(A)$, so we can do this in parallel for all states without conflict. 
For $p\in\{4,\ldots,21\}$ the iteration seems to converge quite rapidly. The reported run for $p=21$ stops after $109$ iterations. At each iteration, the code computes the minimum ratio over states of positive weight, rounded down to six decimal places:
\[
 \alpha=\frac{1}{10^6}\min_{\substack{A\in R\\D(A)>0}}
 \left\lfloor\frac{10^6 T_D(A)}{D(A)}\right\rfloor.
\]
For $p=21$, we get $\alpha=1300544/1000000=1.300544$ after $109$ iterations, which is enough to prove Theorem~\ref{thm:main}.

Although the general algorithm is quite simple, the implementation is slightly more tricky as it involves a lot of details to make it efficient.
First, we perform all computations with integers, rounding in the appropriate direction, so that the inequalities $T_D(A)\ge\alpha D(A)$ hold exactly for the final integer weights. 
We tried to encode everything as tightly as possible to minimize memory consumption as it is the main limiting factor. 
In particular, storing a state requires only $64\times3=192$ bits (slightly less in fact, but we have to work with multiples of 64 bits anyway). The code also contains a simple implementation of a hash table (with open addressing and linear probing) to store the states, their weights and the next weight ($D(A)$ and $T_D(A)$), which allows us to save even more space and improve cache locality. 
We also use specific encodings of the states allowing to use bitwise operations for some of the manipulations that would have required using more costly loops otherwise.

The reported run for $p=21$ has the following size.
\begin{center}
\begin{tabular}{@{}lr@{}}
\toprule
Quantity&Reported value\\
\midrule
States in $R$&$405\,488\,434$\\
States with positive integer weight&$394\,977\,967$\\
Iterations&$109$\\
Memory use&$20.0$ GiB\\
Running time\footnotemark&18 h 24 min\\
\bottomrule
\end{tabular}
\end{center}
\footnotetext{This is the running time with $40$ threads on a machine with $250$ cores, so one might expect it to be roughly $40$ times slower on a single-core machine.}
The run reports that $D(A_0)>0$, and its final ratio computation establishes 
\begin{equation}\label{eq:integer-check}
T_D(A)\ge\frac{1300544}{1000000}D(A)\qquad(A\in R).
\end{equation}
For states of zero weight, these inequalities follow from nonnegativity. These are the computational claims used in Proposition~\ref{prop:computation}.

\section{Open problems}\label{sec:open}

\paragraph{The number of square-free words.}
In the proof of Theorem~\ref{thm:main}, we could have taken $\beta=1.2685$, which improves slightly on the constant $5/4$.
We conjecture that the worst case is when all the lists are identical: for every sequence $L$ of lists of size three, the number of square-free words of length $n$ compatible with $L$ would be at least the number of ternary square-free words of length $n$. 
In particular, the constant $5/4$ could then be replaced by the growth rate $1.3017\ldots$ of ternary square-free words~\cite{shurgrowthrate}.
Even in the online setting of Theorem~\ref{thm:onlineLists}, it seems that the strategy of the adversary that minimizes the number of square-free words of length $n$ that we can construct is to always give the same list.

\paragraph{The erase-repetition game.}
The technique of~\cite{R} was also used in~\cite{ROSENFELDAnnGame} to show that Ann wins the nonrepetitive game over four letters and the erase-repetition game over six letters. Whether Ann wins the erase-repetition game over five or four letters is still open. 
Experiments done in~\cite{ROSENFELDAnnGame} suggest that the technique could prove that she wins over five letters, the limit being again mostly computational, so improvements similar to the one presented here might settle this question.

\paragraph{List packing.}
Like the Thue choice number, the next problem comes from graph colouring. The list packing number $\chi^\star_\ell(G)$ of a graph $G$ is the least $k$ such that, for every assignment of lists of size $k$ to the vertices of $G$, there are $k$ proper colourings from these lists that differ at every vertex~\cite{packingListColourings}.

Inspired by this notion, Grytczuk proposed the following strengthening of the Thue choice number of the path. Let $\pi^\star_{\mathrm{ch}}(P)$ be the least integer $k$ such that, for every $k\times n$ matrix whose columns consist of $k$ distinct letters, one can permute the letters inside each column so that every row is a square-free word. 
By compactness, this amounts to asking that, for every sequence $L$ of lists of size $k$, there are $k$ infinite square-free words respecting $L$ that differ at every position. 
Each row is a square-free word respecting the lists given by the columns, so $\pi_{\mathrm{ch}}(P)\le\pi^\star_{\mathrm{ch}}(P)$. Grytczuk conjectures that $\pi^\star_{\mathrm{ch}}(P)=3$, which would strengthen the existence part of Theorem~\ref{thm:main}; the case $k=3$ was asked at the open problem session of WORDS~2025, as filling a $3\times n$ chessboard~\cite{words2025open}. 
He also informed me that $\pi^\star_{\mathrm{ch}}(P)\le6$ (personal communication).

\section*{Acknowledgements and use of AI}
The main technique behind this article was developed in~\cite{R}, with the
objective of settling this problem, and was later used
in~\cite{cubeFreeWithVuong}. It would have required computational power that
is currently out of reach. Since~\cite{R}, I spent a considerable amount of
time trying to do what is done in the present article, that is, to find a way
to merge the states of the construction of~\cite{R} that greatly reduces
their number, while losing so little information that the growth
constant~$\alpha$ is preserved. I gave Claude (Anthropic) and ChatGPT
(OpenAI) the same prompts: the goal of pushing the computation of~\cite{R} to
an unbounded alphabet, followed by the suggestion of merging similar states
while keeping their worst behaviour.

After a couple of exchanges, Claude found such a way of merging the states:
for each of the last $p$ positions, remember only whether it carries one of
the few most recent letters (this became the last two letters in the article),
together with the length of the repetition in progress for each period
$q\le p$. Claude also gave the argument showing that this abstraction is
sound. ChatGPT, on the other hand, did not produce an answer: the
conversation stopped several times. In both cases, a single interaction used
up the whole token budget, but with different consequences.

Claude then wrote a first, unreadable version of the code that carried out
the computations needed to complete the proof\footnote{Because of its
computational limitations, Claude only ran this code for $p\le15$, and
predicted that the trends would continue for larger $p$, with about
$5\cdot10^8$ states at $p=21$. I confirmed this by running the code on a
large machine.}, as well as first drafts of the proofs, which I then
rewrote. The states used in the present article are considerably simpler
than the ones found by Claude, but they are based on the same idea. I then
rewrote the code in a more readable and verifiable way, with assistance from
Claude. Note that the original unreadable code of Claude is slightly faster than the code provided with this article (it does not compute exactly the same thing, but the main explanation is that it uses optimizations that make the code much more difficult to read and verify for a small performance gain). Both Claude and ChatGPT were used for proofreading, and suggested improvements to the text. ChatGPT was also used to generate its own version of the code to verify that both versions produce the same output. I take full responsibility for the content of this article.
\bibliographystyle{plain} 
\bibliography{biblio}

\end{document}